\documentclass[11pt]{article}

\usepackage{amssymb}
\usepackage{amsthm}
\newtheorem{theorem}{Theorem}

\newtheorem{remark}{Remark}

\newtheorem{conjecture}{Conjecture}
\def\Frac#1#2{\frac{\displaystyle{#1}}{\displaystyle{#2}}}

\usepackage{enumitem}

\begin{document}


\begin{center}
{\Large
{\rm Comments on the paper ``Universal bounds and monotonicity properties of ratios of Hermite and 
Parabolic Cylinder functions"}}
\vspace*{0.5cm}

{\large
{\rm Javier Segura}}
\vspace*{0.5cm}

Departamento de Matem\'aticas, Estad\'{\i}stica y Computaci\'on.

Facultad de Ciencias. 

Universidad de Cantabria. 

39005-Santander, SPAIN. 

javier.segura@unican.es

\vspace*{0.5cm}

------------------------------------------------------------------------------------------------

{\bf Abstract}
\end{center}

In the abstract of \cite{Koch:2020:UBA} we read: ``We obtain so far unproved properties of a ratio involving a class
of Hermite and parabolic cylinder functions." However, we explain how some of the main results in that paper
were already proved in 
\cite{Segura:2012:OBF}, namely the `universal bounds'. 
An error in reference \cite{Segura:2012:OBF} was
discussed in \cite{Koch:2020:UBA} which does not affect the proof given there for 
those `universal bounds'; we fix this erratum easily. We end this note proposing a conjecture regarding the best possible 
upper bound for a certain ratio of parabolic cylinder functions.

\begin{center}
-------------------------------------------------------------------------------------------------
\end{center}





\section{Statement of the results}

We use the notation $U(n,x)$ for parabolic cylinder functions, solutions of the differential equation $y''(x)-(x^2 /4+n) y(x)=0$, 
using the notation of \cite{Temme:2010:PCF}. 

In \cite{Koch:2020:UBA} the following bounds for parabolic cylinder functions 
are obtained (Proposition 3.3), which are presented as new, and which are 
also stated in terms of Hermite functions in Corollary 3.2: 
$$
1<\Frac{D_{\nu-1}(x)^2}{D_\nu (x) D_{\nu -2}(x)}<\Frac{\nu-1}{\nu},\,\nu<0,\,x\in {\mathbb R},
$$

\vspace*{0.3cm}

\footnoterule
The author acknowledges support from Ministerio de Ciencia e Innovaci\'on, project PGC2018-098279-B-I00
(MCIU/AEI/FEDER, UE).

\newpage
where $D_{\nu}(x)=U(-\nu-1/2,x)$. This is the same as saying that

\begin{equation}
\label{eso}
1<F_n(x)=\Frac{U(n,x)^2}{U(n-1,x)U(n+1,x)}<\Frac{n+1/2}{n-1/2},\,n>1/2,\,x\in {\mathbb R}.
\end{equation}

These inequalities are rediscovered in \cite{Koch:2020:UBA} as a consequence of the 
monotonicity of $F_n(x)$, which is the new result in that paper. However the inequalities (\ref{eso}) were already 
proved in \cite{Segura:2012:OBF}, Theorem 11. In this same theorem of \cite{Segura:2012:OBF}
an error was made in a different inequality when declaring its range if validity, but this 
erratum does not affect the proof of statement (\ref{eso}). Surprisingly, the fact 
that the inequalities (\ref{eso}) were already proved 
in \cite{Segura:2012:OBF} was not acknowledged in \cite{Koch:2020:UBA} while 
the erratum in the same theorem of \cite{Segura:2012:OBF} was analyzed in detail, 
even providing graphical information.

We repair the error in Theorem 11 of \cite{Segura:2012:OBF} as follows (but the results on the `universal bounds' are kept intact and its proof remains the same). For convenience we denote 
$$g_{\alpha,\beta}(x)=\Frac{x+\sqrt{4(n-\alpha)+x^2}}{x+\sqrt{4(n-\beta)+x^2}},$$
then:
\begin{theorem}
\label{principal}
Let $F_n(x)=\Frac{U(n,x)^2}{U(n-1,x)U(n+1,x)}$ then the following holds for all real $x$ and for $n>1/2$,  
 except for the first inequality which only holds for $n>3/2$:
$$
\Frac{n-3/2}{n+1/2}g_{-\frac12,\frac32}(x)<
\Frac{n-1/2}{n+1/2}F_n(x)<1<F_n(x)<g_{-\frac32,\frac12}(x)
$$
\end{theorem}

\begin{remark}
Because the function $g_{-\frac32,\frac12}(x)$ is monotonically decreasing in ${\mathbb R}$ we have 
that $F_n(x)<g_{-\frac32,\frac12}(x_0)$ if $x>x_0$ and, in particular 
$$F_n(x)<\displaystyle\sqrt{\Frac{n+\frac32}{n-\frac12}},\, x>0.$$ 

This inequality appeared in Theorem 11 of \cite{Segura:2012:OBF} as being valid for all $x$, when 
it must be restricted to $x>0$. This is the error discussed in \cite{Koch:2020:UBA}. 
However, the so-called `universal bounds' (\ref{eso}) remain intact for $n>1/2$ and all real $x$.
\end{remark}

\begin{remark}
The first inequality in Theorem 11 of \cite{Segura:2012:OBF} is weaker that $F_n(x)>1$, $n>1/2$, and so it must be
 omitted. Instead we include a new inequality 
(which is sharper for $x<0$ than for $x>0$).
\end{remark}

\section{Detailed proof of the theorem}

Having clarified which was the erratum analyzed in detail in \cite{Koch:2020:UBA}, all that remains, for completeness, 
is to recall the proof of the theorem. We provide few more details than 
in \cite{Segura:2012:OBF} in order to
settle the issue completely. We stress again that the `universal bounds' 
and the proof for these results remains the same.

The starting point are the following two results, which are proven in \cite{Segura:2012:OBF} (Theorems 9 and 10) using properties of
the Riccati equations satisfied by the ratios of parabolic cylinder functions, supplemented with the use of the three term
recurrence relation. In the following, we denote 
$$
h_n(x)=\Frac{U(n,x)}{U(n-1,x)}.
$$

\begin{theorem}
\label{bo1}
For $n>1/2$ and $x>0$ the following holds
$$
\Frac{2}{x+\sqrt{4n+2+x^2}}<h_n(x)<\Frac{2}{x+\sqrt{4n-2+x^2}}.
$$
\end{theorem}

The upper bound is one of the characteristic roots of the Riccati equation satisfied by the ratio $h_n(x)$ while the lower
bound is obtained from the upper bound and applying the backward recurrence relation $h_n(x)=(x+(n+1/2)h_{n+1}(x))^{-1}$: using
the upper bound we have that $h_n(x)^{-1}=x+(n+1/2)h_{n+1}(x)<(x+\sqrt{4n+2+x^2})/2$ if $n>-1/2$ (notice that the upper bound
of Thm. \ref{bo1} is valid for $n>1/2$ but that we have the shift $n\rightarrow n+1$). Therefore, if $h_n(x)>0$ (which holds for all $x$ if $n>1/2$ but also for $n\in (-1/2,1/2]$ if $x>0$) we obtain the lower bound.

\begin{remark}
A simple analysis of the Riccati equation suffices to check that $h_n(x)>0$ for all $x$ if $n>1/2$, but also for $n\in(-1/2,1/2]$ when
 $x>0$. Therefore, as mentioned, the first inequality in Theorem \ref{bo1} also holds for $n>-1/2$ when $x>0$, and it turns to an equality when 
$n=-1/2$.
\end{remark}

\begin{theorem}
\label{bo2}
For $n>3/2$ and $x>0$ the following holds
$$
\Frac{x+\sqrt{4n-6+x^2}}{2n-1}<h_n(-x)<\Frac{x+\sqrt{4n-2+x^2}}{2n-1}.
$$
The upper bound also holds for $n\in (1/2,3/2]$.
\end{theorem}

Again, the upper bound is obtained from the analysis of the Riccati equation and the lower bound is obtained from
 the recurrence relation, but now applied in the forward direction, that is: $h_n(x)=(-x+h_{n-1}(x)^{-1})/(n-1/2)$
 (which explains why the lower bound has a smaller
range of validity).

As explained in Remark 4 of \cite{Segura:2012:OBF}, in fact, the previous two
 theorems hold for all real $x$. Indeed, we observe that if in the upper bound of Theorem \ref{bo2} we replace $x$ by $-x$ we obtain the 
 expression of the upper bound of Theorem \ref{bo1}; this means that the upper bound of Theorem \ref{bo1} holds for all real $x$.
Then, applying the recurrence relation the lower bound of Theorem \ref{bo1} also holds for all real $x$ and $n>1/2$. 
Similarly for the lower bound obtained in Theorem \ref{bo2} by applying forward recurrence.
Therefore, we can write the following:

\begin{theorem}
\label{prev}
For all real $x$ the following holds:
\begin{enumerate}
\item{}$
h_n(x)<\Frac{2}{x+\sqrt{4n-2+x^2}},\,n>\frac12
$
\item{}$
h_n(x)>\Frac{2}{x+\sqrt{4n+2+x^2}},\,n>\frac12 \mbox{ (also true for $n>-1/2$ when $x>0$)}
$
\item{}$
h_n(x)>\Frac{n-3/2}{n-1/2}\Frac{2}{x+\sqrt{4n-6+x^2}},\,n>\frac32
$
\end{enumerate}
\end{theorem}

Note that the last inequality in this theorem is the lower bound in Thm. \ref{bo2} but with $x$ replaced by $-x$.

From this last theorem, Theorem \ref{principal} follows using that $F_n(x)=h_{n}(x)/h_{n+1}(x)$:

\begin{enumerate}[label=(\alph*)]
\item{}Considering (1) of Theorem \ref{prev} and (2) but with $n$ replaced by $n+1$ we have $F_n(x)<g_{-\frac32,\frac12}(x)$.
\item{}Using (1) with $n$ replaced by $n+1$ and (2) we have $F_n(x)>1$, and this
also holds for $n\in (-1/2,1/2]$ when $x>0$.
\item{}Considering (1) and (3) but with $n$ replaced by $n+1$ we have $F_n(x)<(n-1/2)/(n+1/2)$.
\item{}Using (1) with $n$ replaced by $n+1$ and (3) we have $F_n(x)>\Frac{n-3/2}{n-1/2}g_{-\frac12,\frac32}(x)$.
\end{enumerate}

We notice that the function $F_n(x)$ has a shape similar to the functions $g_{\alpha,\beta}(x)$, $\alpha<\beta$ and that the selection 
$\alpha=-1/2$, $\beta=1/2$ gives the exact limit values as $x\rightarrow \pm \infty$. We propose the following conjecture

\begin{conjecture}
$F_n(x)<g_{-\frac12,\frac12}(x)$ for all real $x$ and $n>1/2$. 
\end{conjecture}

\begin{remark}
Numerical evidence suggests that this is indeed an upper bound. If so, the bound $g_{-\frac12,\frac12}(x)$ would be
 the best possible upper bound of the form $g_{\alpha,\beta}(x)$, because the limits as $x\rightarrow \pm \infty$ are the exact values
of $F_n(\pm \infty)$. The bound is also sharp as $n\rightarrow +\infty$.
\end{remark}

\bibliographystyle{amsplain}
\bibliography{pat}

\providecommand{\bysame}{\leavevmode\hbox to3em{\hrulefill}\thinspace}
\providecommand{\MR}{\relax\ifhmode\unskip\space\fi MR }
\providecommand{\MRhref}[2]{%
  \href{http://www.ams.org/mathscinet-getitem?mr=#1}{#2}
}
\providecommand{\href}[2]{#2}
\begin{thebibliography}{1}

\bibitem{Koch:2020:UBA}
Torben Koch, \emph{Universal bounds and monotonicity properties of ratios of
  {H}ermite and parabolic cylinder functions}, Proc. Amer. Math. Soc.
  \textbf{148} (2020), no.~5, 2149--2155. \MR{4078099}

\bibitem{Segura:2012:OBF}
Javier Segura, \emph{On bounds for solutions of monotonic first order
  difference-differential systems}, J. Inequal. Appl. (2012), 2012:65, 17.
  \MR{2915630}

\bibitem{Temme:2010:PCF}
N.~M. Temme, \emph{Parabolic cylinder functions}, N{IST} handbook of
  mathematical functions, U.S. Dept. Commerce, Washington, DC, 2010,
  pp.~303--319. \MR{2655352}

\end{thebibliography}

\end{document}